\documentclass[12pt]{article}
\usepackage{latexsym, amssymb, amsmath, amscd, amsfonts, epsfig, graphicx, colordvi,verbatim,ifpdf}
\usepackage{amsfonts, amsmath, amssymb,extarrows}
\usepackage{amssymb,amsfonts,amsmath,latexsym,epsfig,cite, psfrag,eepic,color}
\usepackage{amscd,graphics}
\usepackage{latexsym, amssymb,  amsmath,amscd, amsfonts, epsfig, graphicx, colordvi,amsthm}

\usepackage{graphicx}
\usepackage{color}
\usepackage{ifpdf}
\usepackage{fancybox}

\usepackage{float}

\newtheorem{thm}{Theorem}[section]

\newtheorem{conj}[thm]{Conjecture}

\newtheorem{lem}[thm]{Lemma}

\def\pf{\noindent{\it Proof.} }
\def\qed{\nopagebreak\hfill{\rule{4pt}{7pt}}
\medbreak}

\numberwithin{equation}{section}

\def\qed{\nopagebreak\hfill{\rule{4pt}{7pt}}
\medbreak}

\newlength{\boxedparwidth}
  {\begin{center} \begin{tabular}{|@{\hspace{.315in}}c@{\hspace{.15in}}|}
                  \hline \\ \begin{minipage}[t]{\boxedparwidth}
                  \setlength{\parindent}{.25in}}%
  {\end{minipage} \\ \\ \hline \end{tabular} \end{center}}

\allowdisplaybreaks
\begin{document}

\begin{center}

 {\large \bf The Unimodality of the Crank on Overpartitions}
\end{center}

\begin{center}
{Wenston J.T. Zang}$^{1}$ and
  {Helen W.J. Zhang}$^{2}$ \vskip 2mm
$^{1}$Institute of Advanced Study of Mathematics\\[2pt]
   Harbin Institute of Technology, Heilongjiang 150001, P.R. China\\[8pt]

   \vskip 2mm

 $^{2}$Center for Applied Mathematics\\[2pt]
Tianjin University, Tianjin 300072, P.R. China\\[5pt]
  $^1$zang@hit.edu.cn, \quad $^2$wenjingzhang@tju.edu.cn
\end{center}

\vskip 6mm \noindent {\bf Abstract.} Let $N(m,n)$ denote the number of partitions of $n$ with rank $m$, and let $M(m,n)$ denote the number of partitions of $n$ with crank $m$. Chan and Mao proved that for any nonnegative integers $m$ and $n$, $N(m,n)\geq N(m+2,n)$ and for any nonnegative integers $m$ and $n$ such that $n\geq12$, $n\neq m+2$, $N(m,n)\geq N(m,n-1)$. Recently, Ji and Zang showed that for $n\geq 44$ and $1\leq m\leq n-1$, $M(m-1,n)\geq M(m,n)$ and for $n\geq 14$ and $0\leq m\leq n-2$, $M(m,n)\geq M(m,n-1)$. In this paper, we analogue the result of Ji and Zang to overpartitions. Note that Bringmann, Lovejoy and Osburn introduced two type of cranks on overpartitions, namely the first residue crank and the second residue crank. Consequently, for the first residue crank $\overline{M}(m,n)$, we show that
$\overline{M}(m-1,n)\geq \overline{M}(m,n)$ for  $m\geq 1$ and $n\geq 3$ and $\overline{M}(m,n)\geq \overline{M}(m,n+1)$ for $m\geq 0$ and $n\geq 1$. For the second residue crank $\overline{M2}(m,n)$, we show that
$\overline{M2}(m-1,n)\geq \overline{M2}(m,n)$ for  $m\geq 1$ and $n\geq 0$ and $\overline{M2}(m,n)\geq \overline{M2}(m,n+1)$ for $m\geq 0$ and $n\geq 1$.  Moreover, let $M_k(m,n)$ denote the number of $k$-colored partitions of $n$ with $k$-crank $m$, which was defined by Fu and Tang. They conjectured that when $k\geq 2$, $M_k(m-1,n)\geq M_k(m,n)$ except for $k=2$ and $n=1$. With the aid of the inequality $\overline{M}(m-1,n)\geq \overline{M}(m,n)$ for $m\geq 1$ and $n\geq 3$, we confirm this conjecture.

 \vskip 0.3cm

 \vskip 0.3cm

\section{Introduction}

The first residue crank and the second residue crank on overpartitions are two important statistics in the theory of partitions.
The main objective of this paper is to investigate the distribution of the first residue and the second residue on overpartitions. Recall that the rank of a partition was defined by Dyson \cite{Dyson-1944} as the largest part minus the number of parts. The crank of a partition was introduced by Andrews and Garvan \cite{Andrews-Garvan-1988} as the largest part if the partition contains no ones,
and otherwise as the number of parts larger than the number of ones minus the number
of ones.

Let $p(n)$ denote the number of partitions of $n$. Dyson \cite{Dyson-1944} first conjectured and later confirmed by Atkin and Swinnerton-Dyer \cite{Atkin-Swinnerton-Dyer-1954} that the rank can provide combinatorial interpretations to the first two Ramanujan's congruences $p(5n+4)\equiv 0\pmod{5}$ and $p(7n+5)\equiv 0\pmod{7}$. Andrews and Garvan \cite{Andrews-Garvan-1988} showed that the crank may combinatorially interpret the Ramanujan's congruences $p(5n+4)\equiv 0\pmod{5}$, $p(7n+5)\equiv 0\pmod{7}$ and $p(11n+6)\equiv 0\pmod{11}$.
Relations between ranks and cranks of partitions have been studied by several authors, for example, Andrews and Lewis \cite{Andrews-Lewis-2000}, Garvan \cite{Garvan-1990}, Lewis \cite{Lewis-1991,Lewis-1992}.

%

Let $N(m,n)$ denote the number of partitions of $n$ with rank $m$, and for $n>1$, let $M(m,n)$ denote the number of partitions of $n$ with crank $m$, and for $n=1$, set
\[M(0,1)=-1,~~~~M(-1,1)=M(1,1)=1,\]
and $M(m,1)=0$ for $m\neq0,-1,1$.

Andrews and Garvan \cite{Andrews-Garvan-1988} derived the following generating function of $M(m,n)$:
\begin{equation}\label{equ-gf-crank}
\sum_{m=-\infty}^\infty \sum_{n=0}^\infty M(m,n)z^mq^n=\frac{(q;q)_\infty}{(zq;q)_\infty(q/z;q)_\infty}.
\end{equation}

Here and throughout the rest of this paper, we adopt the common $q$-series notation \cite{Andrews-1976}:
\begin{align*}
(a;q)_\infty=\prod_{n=0}^\infty(1-aq^n) \quad \text{and} \quad
(a;q)_n&=\frac{(a;q)_\infty}{(aq^n;q)_\infty}.
\end{align*}

Chan and Mao \cite{Chan-Mao-2014} discovered the inequality on $N(m,n)$ as stated below.

\begin{thm}[Chan and Mao]
For any nonnegative integers $m$ and $n$,
\[N(m,n)\geq N(m+2,n),\]
and for any nonnegative integers $m$ and $n$ such that $n\geq12$, $n\neq m+2$,
\[N(m,n)\geq N(m,n-1).\]
\end{thm}

Recently, Ji and Zang \cite{Ji-Zang} proved the following unimodal property of the crank on ordinary partitions.

\begin{thm}[Ji and Zang]\label{main-thm-J-Z}
For $n\geq 44$ and $1\leq m \leq n-1$, we have
\begin{equation}\label{main-ine}
M(m-1,n)\geq M(m,n).
\end{equation}
\end{thm}

Ji and Zang also give the following monotonicity property of $M(m,n)$.

\begin{thm}[Ji and Zang]\label{main-un-2-n}
For any $n\geq 14$ and $0\leq m\leq n-2$, we have
\begin{equation}
M(m,n)\geq M(m,n-1).
\end{equation}
\end{thm}

Our interest in this paper is to consider an analogue of Theorem \ref{main-thm-J-Z} and Theorem \ref{main-un-2-n} for overpartitions. Specifically, we will investigate the unimodal property for the first residue crank and the second residue crank. Recall that Corteel and Lovejoy  \cite{Corteel-Lovejoy-2004} defined an overpartition of $n$ as a partition of $n$ in
which the first occurrence of a part may be overlined.
For example, there are $14$ overpartitions of $4$:
\[\begin{array}{lllllllllllll}
(4)&(\overline{4})&(3,1)&(\overline{3},1)&(3,\overline{1})&(\overline{3},\overline{1})&(2,2)\\[5pt]
(\overline{2},2)&(2,1,1)&(\overline{2},1,1)&(2,\overline{1},1)&(\overline{2},\overline{1},1)& (1,1,1,1)&(\overline{1},1,1,1)
\end{array}
\]

Analogous to the  crank of an ordinary partition,
Bringmann, Lovejoy and Osburn \cite{Bringmann-Lovejoy-Osburn-2009} defined the first and second residual crank of an overpartition.
The first residual crank of an overpartition is
defined as the crank of the subpartition consisting of  non-overlined parts. The second
residual crank is defined as the crank of the subpartition consisting of all of the even non-overlined
parts divided by two.

For example, for   $\lambda=(\overline{9},9,7,\overline{6},5,5,\overline{4},4,3,\overline{1},1,1)$, the partition consisting of non-overlined parts of $\lambda$ is $(9,7,5,5,4,3,1,1)$.
The first residual crank of $\lambda$ is $4$. The partition formed by even non-overlined parts of $\lambda$ divided by two is $(2)$. So the second residual crank of $\lambda$ is $2$.

Since then, the researches on overpartitions have been extensively studied (see, for example, Corteel \cite{Corteel-2003}, Corteel and Lovejoy \cite{Corteel-Lovejoy-2002,Corteel-Lovejoy-2004}).
In addition, there have been a great wealth of further results regarding the relations between the rank and crank of overpartitions (see, for example, Andrews, Chan, Kim and Osburn \cite{Andrews-Chan-Kim-Osburn-2016}, Jennings-Shaffer \cite{Jennings-2016}, Lovejoy and Osburn \cite{Lovejoy-Osburn-2008,Lovejoy-Osburn-2010}),
their connections to Maass form and mock theta functions (see, for example, Andrews, Dixit, Schultz and Yee \cite{Andrews-Dixit-Schultz-Yee-2016}, Bringmann and Lovejoy \cite{Bringmann-Lovejoy-2007}, Jennings-Shaffer \cite{Jennings-2016}),
and the relations between $\mathrm{spt}$-function and overpartitions (see, for example, Garvan and Jennings-Shaffer \cite{Garvan-Jennings-2014}, Jennings-Shaffer \cite{Jennings-2015-1,Jennings-2015-2}).

Let $\overline{M}(m,n)$ (resp. $\overline{M2}(m,n)$) denote the number of overpartitions of $n$ with first (resp. second)
residual crank equal to $m$. Here we make the appropriate modifications based on the fact that for ordinary
partitions we have ${M}(0,1)=-1$ and ${M}(-1,1) ={M}(1,1) = 1$. For example, the overpartition $ (\overline{7},\overline{5},\overline{2},1)$ contributes a $-1$ to the count of $\overline{M}(0, 15)$ and a $+1$ to both $\overline{M}(-1, 15)$ and $\overline{M}(1, 15)$. For another example, $(\overline{10}$, $9$, $9$, $\overline{7}$, $7$, $\overline{6}$, $5$, $3$, $3$, $\overline{2}$, $2)$ contributes a $-1$ to the count of $\overline{M2}(0,57)$ and a $+1$ to both $\overline{M2}(-1, 57)$ and $\overline{M2}(1, 57)$.

Bringmann, Lovejoy and Osburn \cite{Bringmann-Lovejoy-Osburn-2009} derived the following generating function of $\overline{M}(m,n)$ and $\overline{M2}(m,n)$:
\begin{align}\label{eq-fg-ov-c}
\sum_{m=-\infty}^\infty \sum_{n=0}^\infty \overline{M}(m,n)z^mq^n&=\frac{(-q;q)_\infty(q;q)_\infty}{(zq;q)_\infty(q/z;q)_\infty}.
\\ \label{eq-fg-ov-c-2}
\sum_{m=-\infty}^\infty \sum_{n=0}^\infty \overline{M2}(m,n)z^mq^n&=\frac{(-q;q)_\infty(q;q)_\infty}{(q;q^2)_\infty
(zq;q)_\infty(q/z;q)_\infty}.
\end{align}

From \eqref{equ-gf-crank}, \eqref{eq-fg-ov-c} and \eqref{eq-fg-ov-c-2}, for fixed $m$,  the generating function of $\overline{M}(m,n)$ is equal to
\begin{align}\label{GF-over-Cr}
\sum_{n=0}^\infty \overline{M}(m,n)q^n=(-q;q)_\infty \sum_{n=0}^\infty {M}(m,n)q^n
\end{align}
and
\begin{align}\label{GF-over-M2Cr}
\sum_{n=0}^\infty \overline{M2}(m,n)q^n=\frac{(-q;q)_\infty}{(q;q^2)_\infty} \sum_{n=0}^\infty {M}(m,n)q^{2n}.
\end{align}

The main result of this paper is  an analogue of Theorem \ref{main-thm-J-Z} for overpartitions as stated below.

\begin{thm}\label{thm-main}
For any $m\geq 1$ and $n\geq 0$,
\begin{equation}\label{ine-main-1}
\overline{M}(m-1,n)\geq \overline{M}(m,n),
\end{equation}
except for $(m,n)=(1,1)$ or $(m,n)=(1,2)$.
\end{thm}

\begin{thm}\label{thm-main-m2}
For any $m\geq 1$ and $n\geq 0$,
\begin{equation}\label{ine-main-1}
\overline{M2}(m-1,n)\geq \overline{M2}(m,n).
\end{equation}
\end{thm}

The proofs of Theorem \ref{thm-main} and Theorem \ref{thm-main-m2} are based on the following Lemma. For the remainder of this paper, let $\{b_n\}_{n=0}^\infty$ and $\{c_n\}_{n=0}^\infty$ be any sequence of nonnegative integers but not necessarily the same sequence in different equations.

\begin{lem}\label{lem-2}
The generating function of $M(m-1,n)-M(m,n)$ can be decomposed as follows.
For $m=1$,
\begin{align}\label{equ-m-0-1-dif-f1-h1}
\sum_{n=0}^\infty (M(0,n)-M(1,n))q^n&=(1-q)^2+q^2(1-q)(1-q^5)(-1+q^2+q^3+q^4-q^5)
\nonumber\\
&\quad\quad+(1-q)\sum_{n=0}^\infty b_nq^n+\sum_{n=0}^\infty c_nq^n.
\end{align}
For $m=2$,
\begin{equation}\label{equ-m2-m1-f2-h2}
\sum_{n=0}^\infty (M(1,n)-M(2,n))q^n=q(1-q)(1-q^3)+(1-q)\sum_{n=0}^\infty b_nq^n+\sum_{n=0}^\infty c_nq^n.
\end{equation}
For any $m\geq 3$,
\begin{equation}\label{equ-m2-m1-fm-hm}
\sum_{n=0}^\infty (M(m-1,n)-M(m,n))q^n=(1-q)\sum_{n=0}^\infty b_nq^n+\sum_{n=0}^\infty c_nq^n.
\end{equation}
\end{lem}

 In light of Lemma \ref{lem-2} and the relations \eqref{GF-over-Cr} and \eqref{GF-over-M2Cr}, we give a proof of Theorem \ref{thm-main} and Theorem \ref{thm-main-m2}. Moreover, we also obtain the following monotonicity property for $\overline{M}(m,n)$ and $\overline{M2}(m,n)$, which is an overpartition analogue of Theorem  \ref{main-un-2-n}.

\begin{thm}\label{thm-main-over-m-m2-n}
For any $m\geq 0$ and $n\geq 1$,
\begin{equation}\label{ine-main-over-mono-1}
\overline{M}(m,n)\geq \overline{M}(m,n-1)
\end{equation}
and
\begin{equation}\label{ine-main-over-mono-2}
\overline{M2}(m,n)\geq \overline{M2}(m,n-1).
\end{equation}
\end{thm}

As an application, we will use Theorem  \ref{thm-main} to prove the unimodality of $k$-crank which was conjectured by Fu and Tang \cite{Fu-Tang-2018}.
Recall that the $k$-colored partition is a $k$-tuple of partitions $\lambda=(\lambda^{(1)},\lambda^{(2)},\ldots,\lambda^{(k)})$. For $k\geq 2$, Fu and Tang \cite{Fu-Tang-2018} defined the $k$-crank of $k$-colored partition as follows:
\begin{equation}
k\text{-crank}(\lambda)=\ell(\lambda^{(1)})-\ell(\lambda^{(2)}),
\end{equation}
where $\ell\left(\pi^{(i)}\right)$ denotes the number of parts in $\pi^{(i)}$.

Let $M_k(m,n)$ denote the number of $k$-colored partitions of $n$ with $k$-crank $m$.
The generating function of $M_k(m,n)$ can be derived by Bringmann and Dousse \cite{Bringmann-Dousse-2016}:
\begin{equation}\label{gen-mkmn}
\sum_{n=0}^\infty \sum_{m=-\infty}^\infty M_k(m,n) z^mq^n=\frac{(q;q)_\infty^{2-k}}{(zq;q)_\infty (z^{-1}q;q)_\infty}.
\end{equation}

Fu and Tang \cite[Conjecture 4.1]{Fu-Tang-2018} raised the following conjecture on $\{M_k(m,n)\}_{m=-n}^n$:

\begin{conj}[Fu and Tang, 2018]\label{conj-f-t}
For $n\geq 0$ and $k\geq 2$, the sequence of $\{M_k(m,n)\}_{m=-n}^n$ is unimodal except for $n=1$, $k=2$.
\end{conj}

Fu and Tang also pointed out that using the asymptotic formula in \cite{Bringmann-Manschot} due to Bringmann and Manschot \cite{Bringmann-Manschot}, one may give an asymptotic proof of Conjecture \ref{conj-f-t}.

In this paper, we confirm Conjecture \ref{conj-f-t} with the aid of Theorem \ref{thm-main}.

This paper is organized as follows: In Section 2, we give a proof of Lemma \ref{lem-2} with the aid of Theorem \ref{main-thm-J-Z}. Section 3 is devoted to prove Theorem \ref{thm-main}. In Section 4,  we demonstrate that Theorem \ref{thm-main-m2} can be deduced by Theorem \ref{thm-main}. In Section 5, we prove Theorem \ref{thm-main-over-m-m2-n}. Finally, in Section 6, Conjecture \ref{conj-f-t} will be confirmed with the aid of Theorem \ref{thm-main}.

\section{The transformation of the generating function of $M(m-1,n)-M(m,n)$}
In this section, we give a proof of  Lemma \ref{lem-2} in light of Theorem \ref{main-thm-J-Z}.

{\noindent \it Proof of Lemma \ref{lem-2}.} From Theorem \ref{main-thm-J-Z}, it is routine to check that
\begin{align}
&\sum_{n=0}^\infty (M(0,n)-M(1,n))q^n\nonumber\\
&=1-2 q+q^3+q^4-q^7-q^9+q^{10}-q^{11}+2 q^{12}-q^{13}\nonumber\\[7pt]
&\quad\quad +2 q^{14}-q^{15}+2 q^{16}-2 q^{17}+3 q^{18}-3 q^{19}+3 q^{20}-2 q^{21}\nonumber\\[7pt]
&\quad\quad+3 q^{22}-3 q^{23}+6 q^{24}-4 q^{25}+6 q^{26}-2 q^{27}+7 q^{28}-4 q^{29}\nonumber\\[7pt]
&\quad\quad+11 q^{30}-5 q^{31}+12 q^{32}-3 q^{33}+13 q^{34}-4 q^{35}+20 q^{36}-6 q^{37}\nonumber\\[7pt]
&\quad\quad+22 q^{38}-q^{39}+27 q^{40}-3 q^{41}+37 q^{42}-q^{43}+\sum_{n=44}^\infty b_nq^n.
\end{align}
Set
\begin{align*}
f(q)&=q^{10}+q^{14}+2q^{16}+3q^{18}+2q^{20}+3q^{22}+4q^{24}+2q^{26}\\
&\quad\quad+4q^{28}
+5q^{30}+3q^{32}+4q^{34}+6q^{36}+q^{38}+3q^{40}+q^{42}
\end{align*}
and
\begin{align*}
h(q)&=q^{14}+q^{20}+2q^{24}+4q^{26}+3q^{28}
+6q^{30}\\
&\quad\quad+9q^{32}+9q^{34}+14q^{36}+21q^{38}+24q^{40}+36q^{42}.
\end{align*}
It is trivial to check that
\begin{align}
\sum_{n=0}^\infty (M(0,n)-M(1,n))q^n&=(1-q)^2+q^2(1-q)(1-q^5)(-1+q^2+q^3+q^4-q^5)\nonumber\\[7pt]
&\quad\quad +(1-q)f(q)+h(q)+\sum_{n= 44}^\infty b_nq^n.
\end{align}
Clearly $f(q)$ and $h(q)$ have nonnegative coefficients. This yields \eqref{equ-m-0-1-dif-f1-h1}.

We next assume that $m=2$. Similar as above, by Theorem \ref{main-thm-J-Z}, it can be checked that
\begin{align}\label{equ-m2-m1-f2-h2-abcd-a}
\sum_{n=0}^\infty (M(1,n)-M(2,n))q^n&=q-q^2-q^4+q^5+q^7+q^9-q^{10}+q^{11}-q^{12}+q^{13}-q^{14}
\nonumber\\[3pt]
&\quad\quad+2 q^{15}-q^{16}+3 q^{17}-q^{18}+4 q^{19}-q^{20}+5 q^{21}-q^{22}\nonumber\\[3pt]
&\quad\quad+6 q^{23}-q^{24}+8 q^{25}-q^{26}+\sum_{n=27}^\infty b_nq^n.
\end{align}
Set
\[f(q)=q^9+q^{11}+q^{13}+q^{15}+q^{17}+q^{19}+q^{21}+q^{23}+q^{25}\]
and
\[h(q)=q^7+q^{15}+2q^{17}+3q^{19}+4q^{21}+5q^{23}+7q^{25}.\]
Together with \eqref{equ-m2-m1-f2-h2-abcd-a}, we see that
\begin{align*}
\sum_{n=0}^\infty (M(1,n)-M(2,n))q^n=q(1-q)(1-q^3)+(1-q)f(q)+h(q)+\sum_{n=27}^\infty b_nq^n.
\end{align*}
This yields \eqref{equ-m2-m1-f2-h2}.

For $3\leq m\leq 7$, it is routine to check that  \eqref{equ-m2-m1-fm-hm} is valid. We next assume that $m\geq 8$. From Theorem \ref{main-thm-J-Z}, we see that for $m\geq 43$ and $n\geq m+1$, we see that $M(m-1,n)-M(m,n)\geq 0$. It can be checked that when $8\leq m\leq 42$ and $44\geq n\geq m+1$, $M(m-1,n)-M(m,n)\geq 0$ also holds. When $n=m$ or $m-1$, we have
\[M(m-1,m-1)-M(m,m-1)=1\]
and
\[M(m-1,m)-M(m,m)=-1.\]
Note that for $n\leq m-2$, $M(m-1,n)=M(m,n)=0$. Thus
\begin{equation}
\sum_{n=0}^\infty \left(M(m-1,n)-M(m,n)\right)q^n=q^{m-1}(1-q)+\sum_{n=m+1}^\infty b_nq^n.
\end{equation}
 This completes the proof.\qed

\section{The proof of Theorem \ref{thm-main}}

This section is devoted to give a proof of Theorem \ref{thm-main}. The main difficulty of the proof is to deal with the case $m=1$. To solve this difficulty, we need two more lemmas.

\begin{lem}\label{lem-3-2}
The coefficients of $q^n$ in
\[(1-q)^2(-q;q)_\infty\]
is nonnegative for all $n\geq 0$ except for $n=1$ or $4$.
\end{lem}

\pf It is clear that
\begin{equation}\label{equ-1-1-q-q-inf-2}
(1-q)^2(-q;q)_\infty=(1-q)(-q;q)_\infty-q(1-q)(-q;q)_\infty.
\end{equation}

We further expand the $(1-q)(-q;q)_\infty$ and $q(1-q)(-q;q)_\infty$ as follows.

We first consider $(1-q)(-q;q)_\infty$. Let $d(n)$ denote the number of distinct partitions of $n$. Clearly
\begin{equation}\label{equ-gf-dis}
\sum_{n=0}^\infty d(n)q^n=(-q;q)_\infty.
\end{equation}
Moreover,  we may classify the set of distinct partitions based on the largest part. To be specific, for $j\geq 1$, let $d_j(n)$ denote the number of distinct partitions of $n$ with the largest part is equal to $j$. Then it is clear that
\begin{equation}\label{equ-gf-dis-dj}
\sum_{n=0}^\infty  d_j(n)q^n=q^j(-q;q)_{j-1}.
\end{equation}
Thus from \eqref{equ-gf-dis} and \eqref{equ-gf-dis-dj} we see that
\begin{equation}
(-q;q)_\infty=\sum_{n=0}^\infty d(n)q^n=\sum_{j=0}^\infty\sum_{n=0}^\infty d_j(n)q^n=
1+\sum_{j=1}^\infty q^j(-q;q)_{j-1}.
\end{equation}
Thus we have
\begin{align}\label{equ-gef-1-q-dis}
&(1-q)(-q;q)_\infty\nonumber\\[3pt]
=&1-q+\sum_{j=1}^\infty q^j(-q;q)_{j-1}-\sum_{j=1}^\infty q^{j+1}(-q;q)_{j-1}\nonumber\\[3pt]
=&1-q+q+\sum_{j=2}^\infty q^j(-q;q)_{j-1}-\sum_{j=2}^\infty q^{j}(-q;q)_{j-2}\nonumber\\[3pt]
=&1+\sum_{j=2}^\infty q^{2j-1}(-q;q)_{j-2}\nonumber\\[3pt]
=&1+q^3+q^5(1+q)+\sum_{j=4}^\infty q^{2j-1}(1+q)(1+q^2)(-q^3;q)_{j-4}.
\end{align}
Moreover,
\begin{align}\label{equ-gef-1-q-dis-1}
&\sum_{j=4}^\infty q^{2j-1}(1+q)(1+q^2)(-q^3;q)_{j-4}\nonumber\\
=&\sum_{j=4}^\infty q^{2j}(-q^2;q)_{j-3}+\sum_{j=4}^\infty q^{2j+1}(-q^3;q)_{j-4}+\sum_{j=4}^\infty q^{2j-1}(-q^3;q)_{j-4}.
\end{align}
Substituting \eqref{equ-gef-1-q-dis-1} into \eqref{equ-gef-1-q-dis}, we derive that
\begin{align}\label{equ-1-1-q-q-inf}
&(1-q)(-q;q)_\infty\nonumber\\[3pt]
=&1+q^3+q^5(1+q)+\sum_{j=4}^\infty q^{2j}(-q^2;q)_{j-3}+\sum_{j=4}^\infty q^{2j+1}(-q^3;q)_{j-4}+\sum_{j=4}^\infty q^{2j-1}(-q^3;q)_{j-4}.
\end{align}

On the other hand, from \eqref{equ-gef-1-q-dis}, we see that
\begin{equation}\label{equ-gef-1-q-dis-3}
q(1-q)(-q;q)_\infty=q+q^4+q^6(1+q)+\sum_{j=4}^\infty q^{2j}(1+q)(1+q^2)(-q^3;q)_{j-4}.
\end{equation}
We deduce that
\begin{align}\label{equ-gef-1-q-dis-2}
&\sum_{j=4}^\infty q^{2j}(1+q)(1+q^2)(-q^3;q)_{j-4}\nonumber\\
&=\sum_{j=4}^\infty q^{2j}(-q^2;q)_{j-3}+\sum_{j=4}^\infty q^{2j+1}(-q^3;q)_{j-4}+\sum_{j=4}^\infty q^{2j+3}(-q^3;q)_{j-4}.
\end{align}
Substituting \eqref{equ-gef-1-q-dis-2} into \eqref{equ-gef-1-q-dis-3}, we have
\begin{align}\label{equ-1-1-q-q-inf-1}
&q(1-q)(-q;q)_\infty\nonumber\\[3pt]
=&q+q^4+q^6(1+q)+
\sum_{j=4}^\infty q^{2j}(-q^2;q)_{j-3}+\sum_{j=4}^\infty q^{2j+1}(-q^3;q)_{j-4}+\sum_{j=4}^\infty q^{2j+3}(-q^3;q)_{j-4}.
\end{align}
Substituting \eqref{equ-1-1-q-q-inf} and \eqref{equ-1-1-q-q-inf-1} into \eqref{equ-1-1-q-q-inf-2}, we deduce that
\begin{align}\label{equ-gef-1-q2-dis}
&(1-q)^2(-q;q)_\infty\nonumber\\[3pt]
&=1-q+q^3-q^4+q^5-q^7+\sum_{j=4}^\infty q^{2j-1}(-q^3;q)_{j-4}-\sum_{j=4}^\infty q^{2j+3}(-q^3;q)_{j-4}\nonumber\\[3pt]
&=1-q+q^3-q^4+q^5+q^9+q^{12}+\sum_{j=6}^\infty q^{2j-1}(-q^3;q)_{j-6}(q^{j-3}+q^{j-2}+q^{2j-5}).
\end{align}
Clearly \eqref{equ-gef-1-q2-dis} has nonnegative coefficients in $q^n$ for all $n\geq 0$ except for $n=1$ or $n=4$. This completes the proof.\qed

\begin{lem}\label{lem-3-3}
The coefficient of $q^n$ in
\[(1-q)(1-q^5)(-1+q^2+q^3+q^4-q^5)(-q;q)_\infty\]
is nonnegative for all $n\geq 1$.
\end{lem}

\pf From the Euler identity \cite{Andrews-1976},
\begin{equation}\label{euler}
(-q;q)_\infty=\frac{1}{(q;q^2)_\infty},
\end{equation}

we deduce that
\begin{align}\label{equ-lem-3-2-main-trans}
&(1-q)(1-q^5)(-1+q^2+q^3+q^4-q^5)(-q;q)_\infty\nonumber\\[3pt]
&=
\frac{(1-q)(1-q^5)(-1+q^2+q^3+q^4-q^5)}{(q;q^2)_\infty}\nonumber\\[3pt]
&=
\frac{(-1+q^2+q^3+q^4-q^5)}{(1-q^3)(q^7;q^2)_\infty}\nonumber\\[3pt]
&=\frac{q^4}{(1-q^3)(q^7;q^2)_\infty}-
\frac{(1-q^2)}{(q^7;q^2)_\infty}.
\end{align}
Let $p_7(n)$ denote the number of partitions of $n$ with each part odd and not less than $7$. It is clear that
\[\sum_{n=0}^\infty p_{7}(n)q^n=\frac{1}{(q^7;q^2)_\infty}.\]
Moreover, we may divide the set of partitions counted by $p_{7}(n)$ into nonintersecting subsets based on the largest part. To be specific, for odd $j\geq 7$, let $p_{7,j}(n)$ denote the number of partitions counted by $p_{7}(n)$ such that the largest part is equal to $j$. It is easy to see that
\begin{equation*}
\sum_{n=0}^\infty  p_{7, j}(n) q^n=\frac{q^j}{(q^7;q^2)_{(j-5)/2}}.
\end{equation*}
Thus we have
\begin{equation}\label{equ-trans-q7-q2-j}
\frac{1}{(q^7;q^2)_\infty}=1+\sum_{j\geq 7\atop j \ odd}  \sum_{n=0}^\infty  p_{7, j}(n) q^n=1+\sum_{j\geq 7\atop j \text{ odd}}\frac{q^j}{(q^7;q^2)_{(j-5)/2}}.
\end{equation}
Using \eqref{equ-trans-q7-q2-j}, we may transform the second term in \eqref{equ-lem-3-2-main-trans} as follows.
\begin{align}\label{equ-lem-3-2-tans-1}
\frac{1-q^2}{(q^7;q^2)_\infty}
&=1+\sum_{j\geq 7\atop j \text{ odd}}\frac{q^j}{(q^7;q^2)_{(j-5)/2}}-q^2-\sum_{j\geq 7\atop j \text{ odd}}\frac{q^{j+2}}{(q^7;q^2)_{(j-5)/2}}\nonumber\\[3pt]
&=1-q^2+\frac{q^7}{1-q^7}+\sum_{j\geq9\atop j \text{ odd}}\frac{q^j}{(q^7;q^2)_{(j-5)/2}}-\sum_{j\geq 9\atop j \text{ odd}}\frac{q^{j}}{(q^7;q^2)_{(j-7)/2}}\nonumber\\[3pt]
&=1-q^2+\frac{q^7}{1-q^7}+\sum_{j\geq 9\atop j \text{ odd}}\frac{q^{j}}{(q^7;q^2)_{(j-7)/2}}\left(\frac{1}{1-q^j}-1\right)\nonumber\\[3pt]
&=1-q^2+\frac{q^7}{1-q^7}+\sum_{j\geq 9\atop j \text{ odd}}\frac{q^{2j}}{(q^7;q^2)_{(j-5)/2}}\nonumber\\[3pt]
&=1-q^2+\frac{q^7}{1-q^7}+\frac{q^{18}}{(1-q^7)(1-q^9)}+\sum_{j\geq 11\atop j \text{ odd}}\frac{q^{2j}}{(q^7;q^2)_{(j-5)/2}}.
\end{align}
We next transform the first term in \eqref{equ-lem-3-2-main-trans} in light of \eqref{equ-trans-q7-q2-j} as follows.
\begin{align}\label{equ-lem-3-2-tans-0}
\frac{q^4}{(1-q^3)(q^7;q^2)_\infty}
&=\frac{q^4}{1-q^3}\left(1+\frac{q^7}{1-q^7}+\frac{q^9}{(1-q^7)(1-q^9)}+\sum_{j\geq 11\atop j\text{ odd}}\frac{q^j}{(q^7;q^2)_{(j-5)/2}}\right)\nonumber\\[3pt]
&=  \frac{q^4}{1-q^3}+\frac{q^{11}}{(1-q^3)(1-q^7)}+\frac{q^{13}}{(1-q^3)(1-q^7)(1-q^9)}\nonumber\\[3pt]
&\quad\quad+
\sum_{j\geq 11\atop j\text{ odd}} \frac{q^{j+4}}{(1-q^3)(q^7;q^2)_{(j-5)/2}}\nonumber\\[3pt]
&=\left(q^4+q^7+\frac{q^{10}}{1-q^3}\right)+\left(1+q^3+\frac{q^6}{1-q^3}\right)
\frac{q^{11}}{1-q^7}\nonumber\\[3pt]
&\quad\quad+\frac{q^{13}}{(1-q^3)(1-q^7)(1-q^9)}+
\sum_{j\geq 11\atop j\text{ odd}} \frac{q^{j+4}}{(1-q^3)(q^7;q^2)_{(j-5)/2}}\nonumber\\[3pt]
&=q^4+q^7+\frac{q^{10}}{1-q^3}+\frac{q^{11}}{1-q^7}+\frac{q^{14}}{1-q^7}
+\frac{q^{17}}{(1-q^3)(1-q^7)}
\nonumber\\[3pt]
&\quad\quad+\frac{q^{13}}{(1-q^3)(1-q^7)(1-q^9)}+
\sum_{j\geq 11\atop j\text{ odd}} \frac{q^{j+4}}{(1-q^3)(q^7;q^2)_{(j-5)/2}}.
\end{align}
Notice that
\begin{equation}\label{equ-lem-3-2-tans-2}
\frac{q^{13}}{(1-q^3)(1-q^7)(1-q^9)}=\frac{q^{13}}{(1-q^7)(1-q^9)}
+\frac{q^{16}}{(1-q^3)(1-q^7)(1-q^9)}.
\end{equation}
Moreover,
\begin{align}\label{equ-lem-3-2-tans-3}
\sum_{j\geq11\atop j\text{ odd}} \frac{q^{j+4}}{(1-q^3)(q^7;q^2)_{(j-5)/2}}&=
\sum_{j\geq11\atop j\text{ odd}} \frac{q^{j+4}}{(1-q^3)(q^7;q^2)_{(j-11)/2}(q^{j-2};q^2)_2}\left(1+
\frac{q^{j-4}}{1-q^{j-4}}\right)\nonumber\\[3pt]
&=\sum_{j\geq11\atop j\text{ odd}} \frac{q^{j+4}}{(1-q^3)(q^7;q^2)_{(j-11)/2}(q^{j-2};q^2)_2}\nonumber\\[3pt]
&\quad\quad+\sum_{j\geq11\atop j\text{ odd}} \frac{q^{2j}}{(1-q^3)(q^7;q^2)_{(j-5)/2}}.
\end{align}
Substituting  \eqref{equ-lem-3-2-tans-2} and \eqref{equ-lem-3-2-tans-3} into \eqref{equ-lem-3-2-tans-0},  we see that
\begin{align}\label{equ-lem-3-2-tans-fin}
\frac{q^4}{(1-q^3)(q^7;q^2)_\infty}&=q^4+\frac{q^{10}}{1-q^3}+\frac{q^{11}}{1-q^7}
+\frac{q^{7}}{1-q^7}
+\frac{q^{17}}{(1-q^3)(1-q^7)}\nonumber\\[3pt]
&\quad\quad+\frac{q^{13}}{(1-q^7)(1-q^9)}
+\frac{q^{16}}{(1-q^3)(1-q^7)(1-q^9)}\nonumber\\[3pt]
&\quad\quad+\sum_{j\geq11\atop j\text{ odd}} \frac{q^{j+4}}{(1-q^3)(q^7;q^2)_{(j-11)/2}(q^{j-2};q^2)_2}\nonumber\\[3pt]
&\quad\quad+\sum_{j\geq11\atop j\text{ odd}} \frac{q^{2j}}{(1-q^3)(q^7;q^2)_{(j-5)/2}}.
\end{align}
By \eqref{equ-lem-3-2-tans-1} and \eqref{equ-lem-3-2-tans-fin}, we see that \eqref{equ-lem-3-2-main-trans} can be rewritten as follows:
\begin{align}\label{equ-lem-3-2-fin-1}
&(1-q)(1-q^5)(-1+q^2+q^3+q^4-q^5)(-q;q)_\infty\nonumber\\[3pt]
&=-1+q^2+q^4+\frac{q^{10}}{1-q^3}
+\frac{q^{17}}{(1-q^3)(1-q^7)}\nonumber\\[3pt]
&\quad\quad+\sum_{j=11\atop j\text{ odd}}^\infty\frac{q^{j+4}}{(1-q^3)(q^7;q^2)_{(j-11)/2}(q^{j-2};q^2)_2}+
\frac{q^{16}}{(1-q^3)(1-q^7)(1-q^9)}\nonumber\\[3pt]
&\quad\quad+\left(\frac{q^{11}}{1-q^7}+\frac{q^{13}}{(1-q^7)(1-q^9)}-
\frac{q^{18}}{(1-q^7)(1-q^9)}\right)\nonumber\\[3pt]
&\quad\quad+\left(\sum_{j=11\atop j\text{ odd}}^\infty\frac{q^{2j}}{(1-q^3)(q^7;q^2)_{(j-5)/2}}-\sum_{j=11\atop j \text{ odd}}^\infty\frac{q^{2j}}{(q^7;q^2)_{(j-5)/2}}\right).
\end{align}
It is clear that
\begin{equation}
\frac{q^{18}}{(1-q^7)(1-q^9)}=\frac{q^{18}}{1-q^7}\left(1+\frac{q^9}{1-q^9}\right)=
\frac{q^{18}}{1-q^7}+\frac{q^{27}}{(1-q^7)(1-q^9)}.
\end{equation}
Thus we have
\begin{align}\label{equ-lem-3-2-fin-2}
&\frac{q^{11}}{1-q^7}+\frac{q^{13}}{(1-q^7)(1-q^9)}-
\frac{q^{18}}{(1-q^7)(1-q^9)}\nonumber\\[3pt]
&=\frac{q^{11}}{1-q^7}-\frac{q^{18}}{1-q^7}
+\frac{q^{13}}{(1-q^7)(1-q^9)}
-\frac{q^{27}}{(1-q^7)(1-q^9)}\nonumber\\[3pt]
&=q^{11}+\frac{q^{13}(1+q^7)}{1-q^9}.
\end{align}
Moreover, it is clear that
\begin{equation}\label{equ-lem-3-2-fin-3}
\sum_{j\geq11\atop j\text{ odd}} \frac{q^{2j}}{(1-q^3)(q^7;q^2)_{(j-5)/2}}-\sum_{j\geq11\atop j \text{ odd}} \frac{q^{2j}}{(q^7;q^2)_{(j-5)/2}}=\sum_{j\geq11\atop j\text{ odd}} \frac{q^{2j+3}}{(1-q^3)(q^7;q^2)_{(j-5)/2}}.
\end{equation}
Substituting \eqref{equ-lem-3-2-fin-2} and \eqref{equ-lem-3-2-fin-3} into \eqref{equ-lem-3-2-fin-1}, we deduce that
\begin{align}\label{equ-lem-3-2-fin}
&(1-q)(1-q^5)(-1+q^2+q^3+q^4-q^5)(-q;q)_\infty\nonumber\\[3pt]
&=-1+q^2+q^4+\frac{q^{10}}{1-q^3}
+\frac{q^{17}}{(1-q^3)(1-q^7)}+\frac{q^{16}}{(1-q^3)(1-q^7)(1-q^9)}\nonumber\\[3pt]
&\quad\quad+\sum_{j\geq11\atop j\text{ odd}} \frac{q^{j+4}}{(1-q^3)(q^7;q^2)_{(j-11)/2}(q^{j-2};q^2)_2}+q^{11}+\frac{q^{13}(1+q^7)}{1-q^9}\nonumber\\[3pt]
&\quad\quad+\sum_{j\geq11\atop j\text{ odd}} \frac{q^{2j+3}}{(1-q^3)(q^7;q^2)_{(j-5)/2}}.
\end{align}
Clearly, \eqref{equ-lem-3-2-fin} has nonnegative coefficients in $q^n$ for all $n\geq 1$. This completes the proof. \qed

We are now in a position to prove Theorem \ref{thm-main}.

\noindent{\it Proof of Theorem \ref{thm-main}.}  From \eqref{GF-over-Cr}, it is clear to see that for fixed integer $m$,
\begin{equation}\label{equ-trans-over-m-m-1-ord}
\sum_{n=0}^\infty \left(\overline{M}(m-1,n)-\overline{M}(m,n)\right) q^n=(-q;q)_\infty \sum_{n=0}^\infty\left( {M}(m-1,n)- {M}(m,n)\right) q^n.
\end{equation}
From \eqref{equ-trans-over-m-m-1-ord}
and Lemma \ref{lem-2}, we see that for $m\geq 3$,
\begin{equation}
\sum_{n=0}^\infty (\overline{M}(m-1,n)-\overline{M}(m,n))q^n=(-q,q)_\infty\left((1-q)\sum_{n=0}^\infty b_nq^n+\sum_{n=0}^\infty c_nq^n\right).
\end{equation}
From \eqref{euler}, we deduce that
\begin{align}\label{equ-m-geq-3-over}
\sum_{n=0}^\infty (\overline{M}(m-1,n)-\overline{M}(m,n))q^n&=\frac{(1-q)\sum_{n=0}^\infty b_nq^n}{(q;q^2)_\infty}
+\frac{\sum_{n=0}^\infty c_nq^n}{(q;q^2)_\infty}\nonumber\\[3pt]
&=\frac{\sum_{n=0}^\infty b_nq^n}{(q^3;q^2)_\infty}
+\frac{\sum_{n=0}^\infty c_nq^n}{(q;q^2)_\infty}.
\end{align}
Clearly, \eqref{equ-m-geq-3-over} has nonnegative coefficients. This yields
\[\overline{M}(m-1,n)\geq \overline{M}(m,n)\]
for $m\geq 3$.

Similar as above, for $m=2$, using Lemma \ref{lem-2}, \eqref{euler} and \eqref{equ-trans-over-m-m-1-ord}, we see that
\begin{align}\label{equ-m-2-over}
\sum_{n=0}^\infty (\overline{M}(1,n)-\overline{M}(2,n))q^n&=(-q,q)_\infty\left(q(1-q)(1-q^3)+(1-q)\sum_{n=0}^\infty b_nq^n+\sum_{n=0}^\infty c_nq^n\right)\nonumber\\[3pt]
&=\frac{q}{(q^5;q^2)_\infty}+\frac{\sum_{n=0}^\infty b_nq^n}{(q^3;q^2)_\infty}
+\frac{\sum_{n=0}^\infty c_nq^n}{(q;q^2)_\infty}.
\end{align}
Clearly \eqref{equ-m-2-over} has nonnegative coefficients. This yields Theorem \ref{thm-main} holds for $m=2$.

Finally, for $m=1$, by  Lemma \ref{lem-2}, \eqref{euler} and \eqref{equ-trans-over-m-m-1-ord}, we deduce that
\begin{align}\label{equ-thm-m-1-final}
&\sum_{n=0}^\infty (\overline{M}(0,n)-\overline{M}(1,n))q^n\nonumber\\
&=q^2(-q;q)_\infty(1-q)(1-q^5)(-1+q^2+q^3+q^4-q^5)\nonumber\\[3pt]
&\quad\quad+(1-q)^2(-q;q)_\infty+\frac{1-q}{(q;q^2)_\infty}\sum_{n=0}^\infty b_nq^n+(-q;q)_\infty \sum_{n=0}^\infty c_nq^n\nonumber\\[3pt]
&=q^2(-q;q)_\infty(1-q)(1-q^5)(-1+q^2+q^3+q^4-q^5)\nonumber\\[3pt]
&\quad\quad+(1-q)^2(-q;q)_\infty+\frac{1}{(q^3;q^2)_\infty}\sum_{n=0}^\infty b_nq^n+(-q;q)_\infty \sum_{n=0}^\infty c_nq^n.
\end{align}
From Lemma \ref{lem-3-2} and Lemma \ref{lem-3-3} we see that
\[q^2(-q;q)_\infty(1-q)(1-q^5)(-1+q^2+q^3+q^4-q^5)+(1-q)^2(-q;q)_\infty\]
has nonnegative coefficients in $q^n$ except for $n=1,2,4$. Thus for $n\neq 1,2,4$, the coefficient of $q^n$ in \eqref{equ-thm-m-1-final} is nonnegative. It is trivial to check that $\overline{M}(0,4)=\overline{M}(1,4)=2$. This completes the proof.\qed

\section{The proof of Theorem \ref{thm-main-m2}}
In this section, we give a proof of Theorem \ref{thm-main-m2} with the aid of Theorem \ref{thm-main}.

\noindent{\it Proof of Theorem \ref{thm-main-m2}.}
By \eqref{GF-over-M2Cr}, it is clear that for any fixed integer $m$,
\begin{align}\label{equ-trans-over-m2-m2-1-ord}
&\sum_{n=0}^\infty \left(\overline{M2}(m-1,n)-\overline{M2}(m,n)\right) q^n\nonumber\\[3pt]
&=\frac{(-q;q)_\infty}{(q;q^2)_\infty} \sum_{n=0}^\infty\left( {M}(m-1,n)- {M}(m,n)\right) q^{2n}\nonumber\\[3pt]
&=\frac{(-q;q^2)_\infty}{(q;q^2)_\infty} \left((-q^2;q^2)_\infty\sum_{n=0}^\infty ({M}(m-1,n)-M(m,n))q^{2n}\right).
\end{align}
From \eqref{equ-trans-over-m-m-1-ord}, we see that
\begin{align}\label{equ-m2-rel-m}
&\sum_{n=0}^\infty (\overline{M2}(m-1,n)-\overline{M2}(m,n))q^n\nonumber\\[3pt]
&=\frac{(-q;q^2)_\infty}{(q;q^2)_\infty} \left(\sum_{n=0}^\infty \left(\overline{M}(m-1,n)-\overline{M}(m,n)\right) q^{2n}\right).
\end{align}
Hence by Theorem \ref{thm-main}, for $m\geq 2$,
\[\sum_{n=0}^\infty (\overline{M2}(m-1,n)-\overline{M2}(m,n))q^n\]
has nonnegative coefficients. This yields Theorem \ref{thm-main-m2} holds for $m\geq 2$.

We next consider the case $m=1$. In this case, by Theorem \ref{thm-main}, it is easy to check that
 \begin{equation}
 \sum_{n=0}^\infty\left(\overline{M}(m-1,n)-\overline{M}(m,n)\right) q^{2n}=1-q^2-q^4+q^6+\sum_{n=8}^\infty b_nq^n.
 \end{equation}
Thus it suffices to show that the coefficients of $q^n$ in
 \[\frac{(-q;q^2)_\infty}{(q;q^2)_\infty}
 (1-q^2-q^4+q^6)\]
is nonnegative for all $n\geq 0$.
Since
\begin{align}\label{equ-m20-asd-sd}
 \frac{(-q;q^2)_\infty}{(q;q^2)_\infty}
 (1-q^2-q^4+q^6)&= \frac{(-q;q^2)_\infty}{(q;q^2)_\infty}
 (1-q^2)(1-q^4)\nonumber\\[3pt]
 &=(-q;q^2)_\infty(1-q^4)\cdot\frac{1+q}{(q^3;q^2)_\infty}.
\end{align}
We claim that
\[(-q;q^2)_\infty(1-q^4)\]
has nonnegative coefficients in $q^n$ for all $n\geq 0$.

Let $sc(n)$ denote the number of distinct odd partitions of $n$. It is clear that
$$(-q;q^2)_\infty=\sum_{n=0}^\infty sc(n)q^n.$$
 We next divide the set of partitions counted by $sc(n)$ into nonintersecting subsets based on the largest part. For $j\geq 1$ odd, let $sc_j(n)$ denote the number of distinct odd partitions of $n$ with largest part equal to $j$. Then
 \begin{equation}
 \sum_{n=0}^\infty sc_j(n)q^n=q^j(-q;q^2)_{(j-1)/2}.
 \end{equation}
So we see that
\begin{equation}\label{equ-sc-dec}
(-q;q^2)_\infty=1+\sum_{j\geq1\atop j\text{ odd}}  q^j(-q;q^2)_{(j-1)/2}.
\end{equation}
Hence we deduce that
\begin{align}\label{equ-sc-1-q4}
(1-q^4)(-q;q^2)
&=1-q^4+\sum_{j\geq1\atop j\text{ odd}}  q^j(-q;q^2)_{(j-1)/2}-\sum_{j\geq1\atop j\text{ odd}}  q^{j+4}(-q;q^2)_{(j-1)/2}\nonumber\\[3pt]
&=1-q^4+q+q^3(1+q)+\sum_{j\geq5\atop j\text{ odd}}  q^j(-q;q^2)_{(j-1)/2}-\sum_{j\geq5\atop j\text{ odd}}  q^{j}(-q;q^2)_{(j-5)/2}\nonumber\\[3pt]
&=1+q+q^3+\sum_{j\geq5\atop j\text{ odd}}  q^j(-q;q^2)_{(j-5)/2}(q^{j-4}+q^{j-2}+q^{2j-6}).
\end{align}
Clearly, \eqref{equ-sc-1-q4} has nonnegative coefficients, and our claim follows.

From the above claim, it is clear to see that \eqref{equ-m20-asd-sd} has nonnegative coefficients of $q^n$ for all $n\geq 0$.
 This completes the proof.\qed

\section{The monotonicity properties of $\overline{M}(m,n)$ and $\overline{M2}(m,n)$}
In this section, we give a proof of Theorem \ref{thm-main-over-m-m2-n}.

\noindent{\it Proof of Theorem \ref{thm-main-over-m-m2-n}.}  For fixed integer $m$,
 \begin{equation}
 \sum_{n=0}^\infty(\overline{M}(m,n+1)-\overline{M}(m,n))q^n=(1-q)\sum_{n=0}^\infty \overline{M}(m,n)q^n.
 \end{equation}
 From \eqref{GF-over-Cr},
  \begin{equation}
 \sum_{n=1}^\infty(\overline{M}(m,n)-\overline{M}(m,n-1))q^n=(1-q)(-q;q)_\infty\sum_{n=0}^\infty {M}(m,n)q^n.
 \end{equation}
 Thus by \eqref{euler}, we have
   \begin{align}
 \sum_{n=1}^\infty(\overline{M}(m,n)-\overline{M}(m,n-1))q^n&=\frac{(1-q)}{(q;q^2)_\infty}\sum_{n=0}^\infty {M}(m,n)q^n\nonumber\\[3pt]
 &=\frac{1}{(q^3;q^2)_\infty}\sum_{n=0}^\infty {M}(m,n)q^n,
 \end{align}
 which clearly has nonnegative coefficients. This yields \eqref{ine-main-over-mono-1}. Using the same argument, we deduce that
 \begin{equation}
 \sum_{n=1}^\infty(\overline{M2}(m,n)-\overline{M2}(m,n-1))q^n=\frac{(-q;q)_\infty}{(q^3;q^2)_\infty} \sum_{n=0}^\infty {M}(m,n)q^{2n},
 \end{equation}
 which implies \eqref{ine-main-over-mono-2}. This completes the proof.\qed

\section{The proof of Conjecture \ref{conj-f-t}}

In this section, we prove Conjecture \ref{conj-f-t} with the aid of Theorem \ref{thm-main}.

\noindent{\it Proof of Conjecture \ref{conj-f-t}.} From \eqref{eq-fg-ov-c} and \eqref{gen-mkmn}, it is clear that when $k\geq 2$
\begin{align}
\sum_{m=-\infty}^\infty\sum_{n=0}^\infty M_k(m,n)z^mq^n&=\frac{(q;q)_\infty^{2-k}}{(zq;q)_\infty(q/z;q)_\infty}\nonumber\\[3pt]
&=\frac{(-q;q)_\infty(q;q)_\infty}{(zq;q)_\infty(q/z;q)_\infty}
\cdot\frac{1}{(q^2;q^2)_\infty(q;q)_\infty^{k-2}}\nonumber\\[3pt]
&=\frac{1}{(q^2;q^2)_\infty(q;q)_\infty^{k-2}}\sum_{m=-\infty}^\infty\sum_{n=0}^\infty \overline{M}(m,n)z^mq^n.
\end{align}
Equating the coefficient of $z^m$ on both sides, we find that for fixed $m$,
\begin{equation}
\sum_{n=0}^\infty M_k(m,n)q^n=
\frac{1}{(q^2;q^2)_\infty(q;q)_\infty^{k-2}}\sum_{n=0}^\infty \overline{M}(m,n)q^n.
\end{equation}
Thus
\begin{equation}\label{equ-diff-mk-om}
\sum_{n=0}^\infty (M_k(m-1,n)-M_k(m,n))q^n=
\frac{1}{(q^2;q^2)_\infty(q;q)_\infty^{k-2}}\sum_{n=0}^\infty (\overline{M}(m-1,n)-\overline{M}(m,n))q^n.
\end{equation}

When $m\geq 2$, by Theorem \ref{thm-main}, we find that $\overline{M}(m-1,n)-\overline{M}(m,n)\geq 0$. Thus by \eqref{equ-diff-mk-om}, we see that $M_k(m-1,n)\geq M_k(m,n)$, as desired.

 We now assume that $m=1$. From Theorem \ref{thm-main} and simple calculation, we have
\begin{equation}\label{eq-exp-om-om-1}
\sum_{n=0}^\infty(\overline{M}(0,n)-\overline{M}(1,n)) q^n=1-q-q^2+q^3+q^5+\sum_{n=6}^\infty b_n q^n.
\end{equation}
 Substituting \eqref{eq-exp-om-om-1} into \eqref{equ-diff-mk-om},
\begin{align}\label{equ-diff-mk-om-2}
&\sum_{n=0}^\infty (M_k(0,n)-M_k(1,n))q^n\nonumber\\[3pt]
=&
\frac{1-q-q^2+q^3+q^5}{(q^2;q^2)_\infty(q;q)_\infty^{k-2}}+\frac{1}
{(q^2;q^2)_\infty(q;q)_\infty^{k-2}}\sum_{n=6}^\infty b_n q^n\nonumber\\[3pt]
=&\left(\frac{1}{(q^4;q^2)_\infty}+\frac{-q+q^3+q^5}{(q^2;q^2)_\infty}\right)
\frac{1}{(q;q)_\infty^{k-2}}+\frac{1}{(q^2;q^2)_\infty(q;q)_\infty^{k-2}}\sum_{n=6}^\infty b_n q^n.
\end{align}
We claim that for $n\geq 2$, the coefficients of $q^n$ in
\[\frac{-q+q^3+q^5}{(q^2;q^2)_\infty}\]
is nonnegative. Recall that Andrews and Merca \cite{Andrews-Merca-2012} proved the following inequality holds for $n>0$:
\begin{equation}
p(n)-p(n-1)-p(n-2)+p(n-5)\leq 0.
\end{equation}
Thus  for $n>0$,
\begin{equation}
p(n)-p(n-1)-p(n-2)\leq 0.
\end{equation}
Therefore
\begin{align}\label{ine-p-n-1-2}
\sum_{n=0}^\infty (p(n-1)+p(n-2)-p(n))q^n=\frac{-1+q+q^2}{(q;q)_\infty}=-1+\sum_{n=1}^\infty b_nq^n.
\end{align}
Setting $q=q^2$ in \eqref{ine-p-n-1-2}, we derive that
\begin{equation}
\frac{-1+q^2+q^4}{(q^2;q^2)_\infty}=-1+\sum_{n=1}^\infty b_nq^{2n}.
\end{equation}
Thus
\begin{equation}\label{equ-diff-mk-om-3}
\frac{-q+q^3+q^5}{(q^2;q^2)_\infty}=-q+\sum_{n=1}^\infty b_nq^{2n+1}.
\end{equation}
This yields our claim.

Substituting \eqref{equ-diff-mk-om-3} into \eqref{equ-diff-mk-om-2}, and notice that the constant term in $1/(q^4;q^2)_\infty$ is equal to $1$,
\begin{equation}\label{equ-diff-mk-om-4}
\sum_{n=0}^\infty (M_k(0,n)-M_k(1,n))q^n=\frac{1-q}{(q;q)_\infty^{k-2}}+\sum_{n=0}^\infty b_nq^n.
\end{equation}
Thus when $k=2$, the coefficients of $q^n$ in \eqref{equ-diff-mk-om-4} is nonnegative for $n\geq 2$. When $k\geq 3$,
\begin{equation}
\sum_{n=0}^\infty (M_k(0,n)-M_k(1,n))q^n=\frac{1}{(q;q)_\infty^{k-3}(q^2;q)_\infty}+\sum_{n=0}^\infty b_nq^n,
\end{equation}
which has nonnegative coefficients. This completes the proof.\qed

 \vskip 0.2cm
\noindent{\bf Acknowledgments.} This work was supported by  the National Science Foundation of China (Grant NO. 11801119).


\begin{thebibliography}{99} \small

\setlength{\itemsep}{-.8mm}
\bibitem{Andrews-1976}
G.E. Andrews, The Theory of Partitions, Addison-Wesley, 1976.

\bibitem{Andrews-Chan-Kim-Osburn-2016}
G.E. Andrews, S.H. Chan, B. Kim and R. Osburn, The first positive rank and crank moments for overpartitions, Ann. Combin. 20 (2) (2016) 193--207.

\bibitem{Andrews-Dixit-Schultz-Yee-2016}
G.E. Andrews, A. Dixit, D. Schultz and A.J. Yee, Overpartitions related to the mock theta function $\omega(q)$, arXiv: 1603.04352.

\bibitem{Andrews-Garvan-1988}
G.E. Andrews and F.G. Garvan, Dyson's crank of a partition, Bull. Amer. Math. Soc. 18 (2) (1988) 167--171.

\bibitem{Andrews-Lewis-2000}
G.E. Andrews and R. Lewis, The ranks and cranks of partitions moduli 2, 3 and 4, J. Number Theory 85 (1) (2000) 74--84.

\bibitem{Andrews-Merca-2012}
G.E. Andrews and M. Merca, The truncated pentagonal number theorem, J. Combin. Theory Ser. A 119 (8) (2012) 1639--1643.

\bibitem{Atkin-Swinnerton-Dyer-1954}
A.O.L. Atkin and P. Swinnerton-Dyer, Some properties of partitions, Proc. Lond. Math. Soc. 66 (1954) 84--106.

\bibitem{Bringmann-Lovejoy-2007}
K. Bringmann and J. Lovejoy, Dyson's rank, overpartitions, and weak Maass forms,
Int. Math. Res. Not. (19) (2007) Art. ID rnm063 34 pp.

\bibitem{Bringmann-Lovejoy-Osburn-2009}
K. Bringmann, J. Lovejoy and R. Osburn, Rank and crank moments for overpartitions, J. Number Theory 129 (7) (2009) 1758--1772.

\bibitem{Chan-Mao-2014}
S.H. Chan and R. Mao, Inequalities for ranks of partitions and the first moment of ranks and cranks of partitions, Adv. Math. 258 (2014) 414--437.

\bibitem{Corteel-2003}
S. Corteel, Particle seas and basic hypergeometric series, Adv. in Appl. Math. 31 (1) (2003) 199--214.

\bibitem{Corteel-Lovejoy-2002}
S. Corteel and J. Lovejoy, Frobenius partitions and the combinatorics of Ramanujan's ${_1\psi}_1$
summation, J. Combin. Theory Ser. A 97 (2002) 179--183.

\bibitem{Corteel-Lovejoy-2004}
S. Corteel and J. Lovejoy, Overpartitions, Trans. Amer. Math. Soc. 356 (2004) 1623--1635.

\bibitem{Bringmann-Dousse-2016}
K. Bringmann and J. Dousse, On Dyson's crank conjecture and the uniform asymptotic behavior of certain inverse theta functions, Trans. Amer. Math. Soc. 368 (5) (2016) 3141--3155.

\bibitem{Bringmann-Manschot}
K. Bringmann and J. Manschot, Asymptotic formulas for coefficients of inverse theta functions,
Commun. Number Theory Phys. 7 (3) (2013) 497--513.

\bibitem{Dyson-1944}
F.J. Dyson, Some guesses in the theory of partitions, Eureka (Cambridge) vol. 8 (1944) 10--15.

\bibitem{Fu-Tang-2018}
S. Fu and D. Tang, On a generalized crank for $k$-colored partitions, J. Number Theory 184 (2018) 485--497.

\bibitem{Garvan-1990}
F.G. Garvan, The crank of partitions mod 8, 9 and 10, Trans. Amer. Math. Soc. 322 (1) (1990) 79--94.

\bibitem{Garvan-Jennings-2014}
F.G. Garvan and C. Jennings-Shaffer, The spt-crank for overpartitions, Acta Arith. 166 (2) (2014) 141--188.

\bibitem{Jennings-2015-1}
C. Jennings-Shaffer, Another SPT crank for the number of smallest parts in overpartitions with even smallest part,
J. Number Theory 148 (2015) 196--203.

\bibitem{Jennings-2015-2}
C. Jennings-Shaffer, Higher order SPT functions for overpartitions, overpartitions with smallest part even, and partitions with smallest part even and without repeated odd parts, J. Number Theory 149 (2015) 285--312.


\bibitem{Jennings-2016}
C. Jennings-Shaffer, Overpartition rank differences modulo 7 by Maass forms, J. Number Theory 163 (2016) 331--358.

\bibitem{Ji-Zang}
K.Q. Ji and W.J.T. Zang, Unimodality of the Andrews-Garvan-Dyson cranks of  partitions, submitted, 	arXiv:1811.07321.



\bibitem{Lewis-1991}
R. Lewis, On the rank and the crank modulo 4, Proc. Amer. Math. Soc. 112 (4) (1991) 925--933.

\bibitem{Lewis-1992}
R. Lewis, On some relations between the rank and the crank, J. Combin. Theory Ser. A 59 (1) (1992) 104--110.

\bibitem{Lovejoy-Osburn-2008}
J. Lovejoy and R. Osburn, Rank differences for overpartitions, Quart. J. Math. 59
(2008) 257--273.

\bibitem{Lovejoy-Osburn-2010}
J. Lovejoy and R. Osburn, $M_2$-rank differences for overpartitions, Acta Arith. 144 (2) (2010) 193--212.

\end{thebibliography}
\end{document}